\documentclass[11pt]{article}
\usepackage{anysize}
\marginsize{1in}{1in}{1in}{1in}
\usepackage{amsfonts,amsmath,amsxtra,euscript, color}
\usepackage{amssymb}
\usepackage{graphicx,longtable,wrapfig}
\usepackage{epstopdf}
\usepackage[english]{babel}

\title{Multivariate Abel--Ruffini}
\date{}

\author{A. Esterov \thanks{National Research University Higher School of Economics \newline Faculty of Mathematics, NRU HSE, 7 Vavilova 117312 Moscow, Russia, {\sc aesterov@hse.ru}. This study (research grant No 14-01-0152) is supported by The National Research University--Higher School of Economics' Academic Fund Program in 2014/2015. Partially supported by RFBR grant 13-01-00755 and the Dynasty Foundation fellowship.} \and G. Gusev \thanks{Moscow Institute of Physics and Technology (State University), Yandex\newline Department of Innovations and High Technology, MIPT,  9 Institutskii per. 141700 Dolgoprudny, Moscow Region, Russia, {\sc gusev@mccme.ru}. Partially supported by RFBR grant 13-01-00755.}}

\newtheorem{thm}{Theorem}

\newtheorem{prop}{Proposition}
\newtheorem{lem}{Lemma}
\newtheorem{conj}{Conjecture}
\newtheorem{exa}{Example}

\newenvironment{dfntn} {\smallskip\noindent{\bf
Definition\/}.}{\smallskip\par}

\newenvironment{exmpl} {\smallskip\noindent{\bf Example\/}.}{\smallskip\par}
\newenvironment{rmrk} {\smallskip\noindent{\bf Remark\/}.}{\smallskip\par}

\newenvironment{prf-1} {\noindent{\em Proof of Lemma~\ref{lem-main}}.}{{ $\Box$}\smallskip\par}
\newenvironment{prf-2} {\noindent{\em Proof of Lemma~\ref{Permutations}}.}{{ $\Box$}\smallskip\par}

\renewcommand{\a}{\alpha }

\renewcommand{\b}{\beta }

\newcommand{\g}{\gamma}
\newcommand{\G}{\Gamma }
\renewcommand{\l}{\lambda }

\newcommand{\p}{\pi}

\newcommand{\MV}{{\rm MV}}

\newcommand{\As}{\mathbf A}
\newcommand{\ff}{\mathbf f}
\newcommand{\xx}{\mathbf x}

\newcommand{\yy}{\mathbf y}

\newcommand{\kk}{\mathbf k}
\newcommand{\aaa}{\mathbf a}
\newcommand{\R}{\mathbb{R}}

\newcommand{\C}{\mathbb{C}}
\newcommand{\N}{\mathbb{N}}

\newcommand{\Z}{\mathbb{Z}}
\newcommand{\ZZ}{\EuScript{Z}}
\renewcommand{\AA}{\As}

\newcommand{\Id}{\mathop{\mathrm{Id}}\nolimits}
\newcommand{\Vol}{\mathop{\mathrm{Vol}}\nolimits}
\newcommand{\conv}{\mathop{\mathrm{conv}}\nolimits}
\newcommand{\CC}{(\mathbb{C}\setminus 0)}

\begin{document}
\maketitle

\begin{abstract}
We generalize the Abel---Ruffini theorem to arbitrary dimension, i.e. classify general square systems of polynomial equations solvable by radicals. In most cases, they reduce to systems whose tuples of Newton polytopes have mixed volume not exceeding~4. The proof is based on topological Galois theory, which ensures non-solvability by any formula involving quadratures and single-valued functions, and the computation of the monodromy group of a general system of equations, which may be of independent interest. \vspace{1ex}

\end{abstract}

\section{Introduction}

The classical Abel---Ruffini theorem states that the general univariate polynomial equation of degree $d$ is solvable by radicals if and only if the number of its solutions is less than five, which is equivalent, by the fundamental theorem of algebra, to $d<5$. For a multivariate polynomial, a natural generalization of the degree is its Newton polytope, and the natural generalization of the fundamental theorem of algebra is the following Kouchnirenko theorem (\cite{bernst}): the square system of general equations with a Newton polytope $A$ has $\Vol A$ solutions, where $\Vol$ is the lattice volume normalized by the volume of the standard simplex.

We describe all lattice polytopes $A$ such that the square system of general equations with the Newton polytope $A$ is solvable by radicals. The answer consists of all polytopes of lattice volume less than 5 (Theorem \ref{Thm_equal_polytops}). In particular, we classify all lattice polytopes of lattice volume not exceeding 4 (Section \ref{Slattice}). Our topological approach ensures that systems with the Newton polytope of volume greater than 4 are not solvable in a much stronger sense: the solution cannot be given by a formula involving quadratures and single-valued functions.

If we do not restrict our attention to systems, all of whose equations have the same Newton polytope, then the classification of solvable systems becomes less straightforward: there do exist solvable systems of general equations with more than 4 solutions. In this generality, we only conjecture the answer (Conjecture \ref{conj1}) and prove it in some special cases (Theorems \ref{Thm_dimension_two}, \ref{mix1}, and \ref{mix2}; in the latter case, we can even prove that the system with $N$ solutions cannot be solved by a formula involving quadratures, single-valued functions and solutions of algebraic equations of degree smaller than $N$).

Our results are based on the computation of the monodromy group for a square system of general equations, or, more generally, the monodromy $\zeta$-function of an arbitrary system (Sections \ref{Smonodr0} and \ref{Smonodr}). We then relate the solvability of the system to the solvability of its monodromy group by means of topological Galois theory. For a polynomial $f\in\C[x],\, f:\C\to\C$, its monodromy group is usually interpreted as the Galois group of the polynomial $f(x)-t\in\C(t)[x]$ (see e.g. \cite{topkhov}), so the univariate version of our approach falls within the scope of the classical Galois theory over function fields. However, for systems of more than two equations, we do not see a natural algebraic interpretation for our topological approach. We implement this approach in two ways: the first one (presented in Section \ref{Smix1}) resembles Ritt's technique and leads to the proof of Theorems \ref{Thm_equal_polytops},  \ref{Thm_dimension_two}, and \ref{mix1}. The second one (presented in Sections \ref{Sproofs2} and \ref{Smix2}) resembles Arnold's and Khovanskii's technique and yields Theorem \ref{mix2}.

\section{Abel---Ruffini theorem for multivariate polynomials}

We first introduce definitions and notation used throughout the paper. Let $\As$ be a tuple of finite sets $A_1,A_2,\ldots,A_n\subset\Z^n$. For every $j=1,2,\ldots,n$, denote a tuple of complex numbers $(c_{j,\aaa},\, \aaa\in A_j)$ by $c_{A_j}$, and the space of all such tuples by $\C^{A_j}$. In this paper, we consider the space $\C^{\As}:=\C^{A_1}\oplus\C^{A_2}\oplus\ldots\oplus\C^{A_n}$ of tuples $c_{\As}=(c_{A_1},c_{A_2},\ldots, c_{A_n})$. This is the configuration space of systems of (Laurent) polynomials $\ff=(f_1,f_2,\ldots, f_n)$, $f_j(\xx)=\sum_{\aaa\in A_j} c_{j,\aaa}\xx^\aaa$, where $\xx^\aaa$ stands for the monomial $x_1^{a_1}\cdot x_2^{a_2}\cdot\ldots\cdot x_n^{a_n}$.

\begin{dfntn}
\begin{enumerate}
\item The  {\it solution of the general system of equations supported at $A_1,A_2\ldots,A_n$,} is the (multivalued) function $F:\C^{\As}\to\CC^n$ whose value at every tuple $c_{\As}\in\C^{\As}$ is the set of solutions for the system of polynomial equations
$$
\sum\nolimits_{\aaa\in A_j} c_{j,\aaa} \xx^\aaa=0,\;\; j=1,2,\ldots,n.\eqno{(*)}
$$

\item The general system is {\it solvable by radicals}, if every point $c_\As\in \C^{\As}$ such that $F(c_\As)$ is finite
admits a Zariski open neighborhood $U$ and a multivalued function $G\colon U\to\CC^n$
on it such that $F(c)\subset G(c)$ for $c\in U$, and $G$ is a composition of rational functions and roots of arbitrary degree. The system is {\it solvable by generalized quadratures}, if $G$ is allowed to be a composition of single-valued functions, roots and taking antiderivatives.
\end{enumerate}
\end{dfntn}

\begin{rmrk} We cannot expect $G$ to be defined on the whole $\C^\As=\C^{A_1}\oplus\C^{A_2}\oplus\ldots\oplus\C^{A_n}$: if $n=1$ and $A_1=\{0,1,2\}$, then the standard formula $\frac{-c_1\pm\sqrt{D}}{2c_2}$ for the roots of the equation $c_2x^2+c_1x+c_0=0$ is defined outside the plane $\{c_2=0\}\subset\C^{\As}$, and there is another formula for the roots defined on this plane outside $\{c_0=0\}$, namely $\frac{2c_0}{-c_1\pm\sqrt{D}}$.
We also cannot expect $F(c)=G(c)$: the Cardano formula for the roots of the general cubic equation gives six values, of which only three are the roots.
\end{rmrk}

\begin{rmrk}\label{monotonicity}
Solvability satisfies the following monotonicity property: if the general system of equations supported at $A_1,A_2,\ldots, A_n\subset \Z^n$ is solvable by radicals or quadratures, then so is the genral system supported at arbitrary subsets $B_1\subset A_1,\, B_2\subset A_2,\,\ldots,B_n\subset A_n$.
\end{rmrk}

A tuple $A_1,A_2,\ldots,A_n\subset\Z^n$ is said to be {\it irreducible}, if $0\in \bigcap_j A_j$ and $\bigcup_j A_j$ generates $\Z^n$.
We restrict our attention to irreducible tuples from here on in the paper, because the general system of equations $\{f_1(\xx) = f_1(\xx)=\ldots = f_n(\xx)=0\}$ supported at an arbitrary tuple $A_1,A_2,\ldots,A_n$ is equivalent to a system supported at an irreducible one. In fact, let $A'_1,A'_2,\ldots,A'_n$ be shifted copies of $A_1,A_2,\ldots,A_n$ that contain $0$. Assume they generate a sublattice $L\subset\Z^n$. If ${\rm rk\, } L<n$, then the system $\{f_j(\xx)=0\}$ is inconsistent. Otherwise, $L$ is the image of the inclusion $\varphi:\Z^n\to \Z^n$ given by a matrix $(\varphi_{ik})$, and the system $\{f_j(\xx)=0\}$ is equivalent to the general system $\{g_j(\yy)=0\}$ supported at the sets $\varphi^{-1}(A'_1),\ldots,\varphi^{-1}(A'_n)$ under the change of variables $y_i=\prod_k x_k^{\varphi_{ki}}$.

\subsection{Equations supported at the same set}

\begin{thm}\label{Thm_equal_polytops}
The general system of equations supported at an irreducible tuple $A_1=A_2=\ldots=A_n=A\subset\Z^n$ is solvable by radicals ($\Leftrightarrow$ in generalized quadratures) if and only if it has at most 4 solutions, i.e. the volume of the convex hull $\conv A$ of the polytope $A$ is at most 4.
Every such set $A$ is contained in one of the following 34 sets or in the set obtained from one of those 34 by applying iteratively
the following procedures: (1) taking the standard cone $B\leadsto \{0,\ldots,0,1\}\;\cup\;B\times\{0\}\;\subset\;\Z^{m+1}$ over $B\in \Z^m$, (2) taking the image under an affine automorphism of the lattice.
\begin{itemize}
\item $n=6,\;\; \Vol(\conv A)=4:$ the circuit $S_6\cup\{(-1,-1,-1,1,1,1)\}$, where $S_n$
is the set of vertices of the standard $n$-dimensional simplex.

\item $n=5,\;\; \Vol(\conv A)=4:$ the circuit $S_5\cup\{(-2,-1,1,1,1)\}$ and the join $(S_1\times S_1)\star(S_1\times S_1)$, where $A\star B$ for $A\subset\Z^m$ and $B\subset\Z^n$ is the union $A\times\{0\}\times\{0\}\cup\{0\}\times B\times \{1\}\subset\Z^m\oplus\Z^n\oplus\Z$.

\item $n=4,\;\; \Vol(\conv A)=4:$
\begin{itemize}
\item the circuits $S_4\cup\{(-2,-1,1,1)\},\;\;
S_4\cup\{(-1,-1,-1,1)\},\;\; S_4\cup\{(-1,-1,-1,2)\},$
\item the prism $S_1\times S_3,$
\item the join $(2S_1)\star(S_1\times S_1),$
\item the sum $(S_1\times S_1)\oplus(S_1\times S_1)$, where $A\oplus B$ for $A\subset\Z^m$ and $B\subset\Z^n$ is the union $A\times\{0\}\cup\{0\}\times B\subset\Z^m\oplus\Z^n$.
\end{itemize}
\item $n=4,\;\; \Vol(\conv A)=3:$ the circuit $S_4\cup\{(-1,-1,1,1)\}$.

\item $n=3,\,\, \Vol(\conv A)=4$ and $3$: the circuits $S_3\cup\{(-1,-1,-1)\},\, S_3\cup\{(1,1,-3)\},\, S_3\cup\{(1,1,-2)\}$, the prism $P=S_2\times S_1$ and the sets $P\cup\{(0,0,2)\},\, \{-1,0,1\}\star\{-1,0,1\},\, D\cup\{(0,0,-1)\},\, D\cup\{(0,0,2)\},\, D\cup\{(1,1,1)\}$ and $D\cup\{(1,1,-1)\}$, where $D$ is the square pyramid $S_3\cup\{(1,1,0)\}$. All of them (except for the first one) are shown below:

\begin{center}
\noindent\includegraphics[width=12cm]{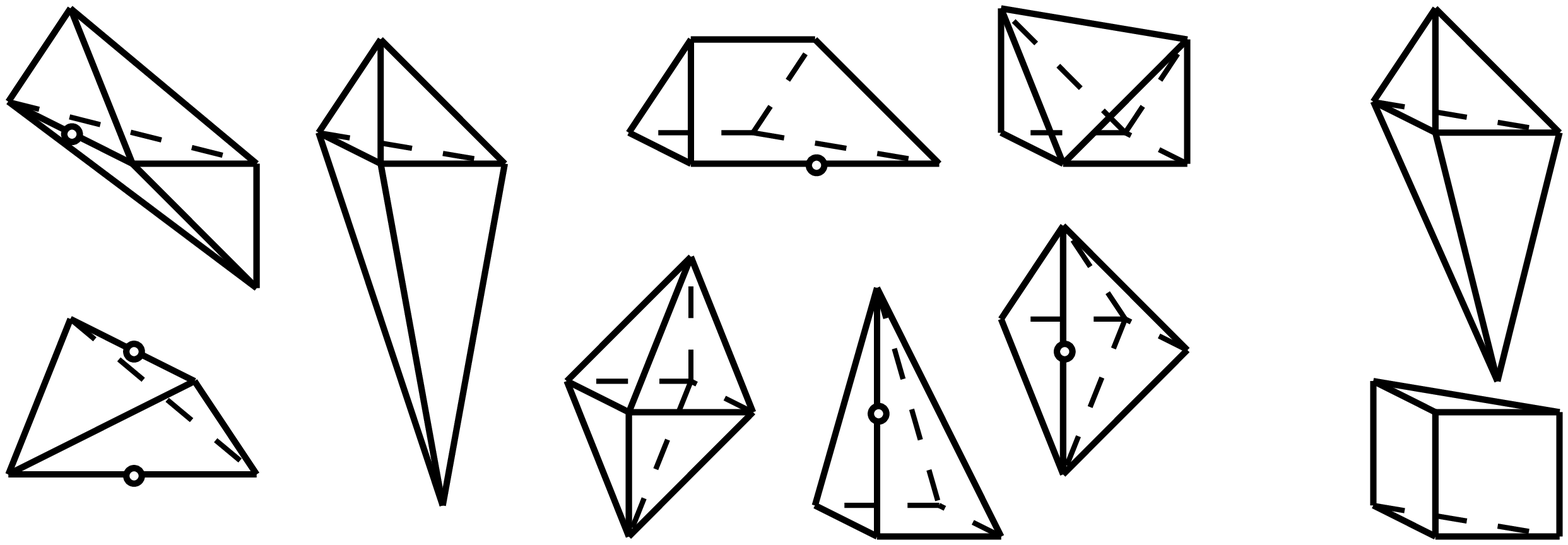}
\end{center}

\item $n=2,\;\; \Vol(\conv A)\leqslant 4$:

\begin{center}
\noindent\includegraphics[width=11cm]{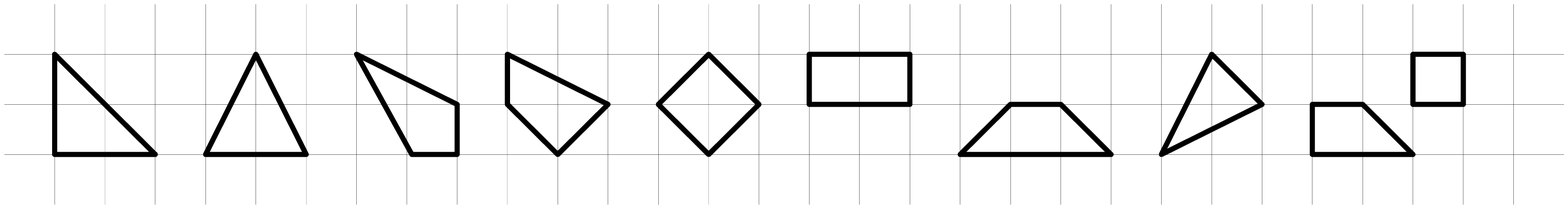}
\end{center}

\item $n=1,\;\; \Vol(\conv A)\leqslant 4: S_1,\, 2S_1,\, 3S_1,\, 4S_1.$ \vspace{1ex}
\end{itemize}
\end{thm}

\begin{exmpl} The roots of the general equation $c_{a_0}x^{a_0}+c_{a_1}x^{a_1}+\ldots+c_{a_n}x^{a_n}=0,\, a_0<a_1<\ldots<a_n,$ can be expressed by radicals in terms of the coefficients $c_{a_0},c_{a_1},\ldots,c_{a_n}$ if and only if $\frac{a_n-a_0}{{\rm GCD}(a_1-a_0,\ldots,a_n-a_0)}\leqslant 4$.
\end{exmpl}

The ``only if'' part of this theorem can be reduced to the subsequent Theorem~\ref{mix1} by means of Example~\ref{lucky1} and monotonicity of solvability. The relation between the number of solutions and the volume is the Kouchnirenko formula. The classification of lattice volume 4 polytopes is proved in Section \ref{Slattice}. The ``if'' part follows from the observation that each of the 34 systems on the list above is solvable (we do not provide 34 obvious explicit formulas here).

\begin{rmrk} The list above only includes irreducible lattice sets, the volume of whose convex hull is at most 4. In order to list all the lattice polytopes of volume at most 4 (up to affine automorphisms of the lattice and taking the standard cone), we should add the empty simplices of volume at most 4 and the pyramides $D_n\oplus(S_1\times S_1),\, D_n\star(S_1\times S_1),\, D_n\oplus(2S_1),\, D_n\star(2S_1)$, where $D_n$, for $n\geqslant 3$, is the (unique) $n$-dimensional volume 2 empty simplex that is not a standard cone over another volume 2 simplex. Recall that an empty, or lattice-free simplex is a simplex that contains no lattice points besides its vertices. Classification of empty simplices of small volume is an interesting question out of the scope of this paper; see e.g. \cite{empty1} for how to classify multidimensional empty simplices.
\end{rmrk}

\subsection{Equations with arbitrary supports}

We now try to drop the assumption $A_1=A_2=\ldots=A_n$.
\begin{prop} Assume that $\Delta\subset\R^n$ is a lattice polytope, and every set $A_i$ in an irreducible tuple $A_1,\ldots,A_n$ equals $\Z^n\cap(d_i\cdot\Delta)$ for some $d_i\in\N$. Then the general system of equations supported at $(A_1,\ldots,A_n)$ is solvable by radicals ($\Leftrightarrow$ in generalized quadratures) if and only if $d_1\cdot\ldots\cdot d_n\cdot\Vol\Delta\leqslant 4$.
\end{prop}
The proof is the same as for Theorem \ref{Thm_equal_polytops}.

\begin{exmpl} The general square system of polynomial equations of degrees $d_1,\ldots,d_n$ is solvable by radicals if and only if $d_1\cdot\ldots\cdot d_n\leqslant 4$.
\end{exmpl}
However, as soon as the convex hulls of the support sets are not homothetic, the question of solvablility becomes drastically more complicated. For instance, an irreducible tuple
$$
A_1=\{0,1,2,3,4\}\times\{0\},\, A_2=\{0\}\times\{0,1,2,3,4\}
$$
gives rise to a general system of equations, which has 16 solutions but is solvable by radicals. This is because this tuple is not reduced in the following sense.

\begin{dfntn} An irreducible tuple $A_1,A_2\ldots,A_n\subset\Z^n$ is said to be {\it reduced}, if the dimension of the convex hull of the union $\cup_{j\in I} A_j$ is greater than $|I|$ for every $I\subsetneq\{1,2,\ldots,n\}$.
\end{dfntn}

From here on in the paper, we restrict our attention to reduced tuples, because the question of solvability of the general system supported at a non-reduced tuple can be reduced to the same question for systems of fewer variables as follows.
For an irreducible tuple of sets $A_1,A_2\ldots,A_n\subset\Z^n$, which is not reduced, a suitable automorphism of the lattice $\varphi:\Z^n\to\Z^n$ sends $A_j$ to $A'_j$ such that $A'_1,A'_2,\ldots,A'_k$ are contained in the first coordinate plane $\Z^k\times\{0,\ldots,0\}\subset\Z^n,\, k<n$. The corresponding change of variables $y_i=\prod_m x_m^{\varphi_{mi}}$, where $(\varphi_{im})$ is the matrix of $\varphi$, sends the general system $f_j(\xx)=0$ supported at $\AA$ to the general system $g_j(\yy)=0$ supported at $\AA':=(A'_1,A'_2,\ldots,A'_k)$, and the latter can be solved in two steps: first, solve the system $g_1(y_1,y_2,\ldots,y_k)=g_2(y_1,y_2,\ldots,y_k)=\ldots=g_k(y_1,y_2\ldots,y_k)=0$, then,  for every solution $(y_1^0,y_2^0,\ldots,y_k^0)$, solve the system $g_j(y^0_1,\ldots,y^0_k,y_{k+1},\ldots,y_n)=0,\, j=k+1,k+2,\ldots,n$. The first system is the general system supported at $A'_1,A'_2,\ldots,A'_k\subset\Z^k$, and the second one is equivalent to the general system supported at the images of $A'_{k+1},A'_{k+2},\ldots,A'_n$ under the projection $\Z^n\to\Z^n/\Z^k\times\{0,0,\ldots,0\}$.

\begin{conj} \label{conj1} The general system of polynomial equations supported at a reduced tuple \linebreak $A_1,A_2,\ldots,A_n$ is solvable by radicals ($\Leftrightarrow$ in generalized quadratures) if and only if it has at most 4 solutions.
\end{conj}

\begin{rmrk}\label{K-B} Recall the Kouchnirenko---Bernstein formula (\cite{bernst}): the number of solutions for the general system of equations supported at a tuple $A_1,A_2,\ldots,A_n$ equals the {\it mixed volume} of the convex hulls of $A_1,A_2,\ldots,A_n$, which is defined as the unique symmetric function of $n$ polytopes in $\R^n$, multilinear with respect to the Minkowski addition $A+B=\{a+b\,|\, a\in A,\, b\in B\}$ and assigning the number $\Vol(A)$ to every tuple of the form $(A,A,\ldots,A)$. Because of this, Conjecture \ref{conj1} raises the question of classification of reduced tuples of polytopes of mixed volume 4. Reduced tuples of mixed volume 1 are classified in \cite{1sol}. Reduced tuples of volume $\leqslant 4$ in $\R^2$ and of volume $\leqslant 2$ in $\R^3$ are classified below.
\end{rmrk}

\begin{rmrk}  The solvability of a general system with at most 4 solutions easily follows from a general argument of elimination theory: let $\{f_j(\xx)=0\}$ be the general system of equations supported at $\AA$, consider $f_j$ as a Laurent polynomial $F_j$ in $x_2,x_3,\ldots,x_n$ with coefficients in $\C[x_1,x^{-1}_1]$, then the mixed resultant of $F_2,F_3,\ldots,F_n$ equals $R\in\C[x_1,x^{-1}_1]$. The roots of $R$ are the first coordinates of the roots of the system $\{f_j(\xx)=0\}$; if there are at most 4 of them, then the equation $R=0$ can be solved by the Ferrari formula. See \cite{1sol} for a more explicit polynomial time algorithm that finds the solution of a general system provided that the solution is unique.
\end{rmrk}

\begin{thm}\label{Thm_dimension_two}
Conjecture~\ref{conj1} is valid for $n=2$. Moreover, the general system of two polynomial equations supported at a reduced pair $A,B\subset\Z^2$ has less than 5 solutions if and only if there exist $G\in SL(\Z^2)$ and $a$ and $b\in\Z^2$, such that the sets $GA+a$ and $GB+b$ are contained in one of the 14 pairs of polygons in Figure 1.

\begin{figure}
\begin{center}
\noindent\includegraphics[width=11cm]{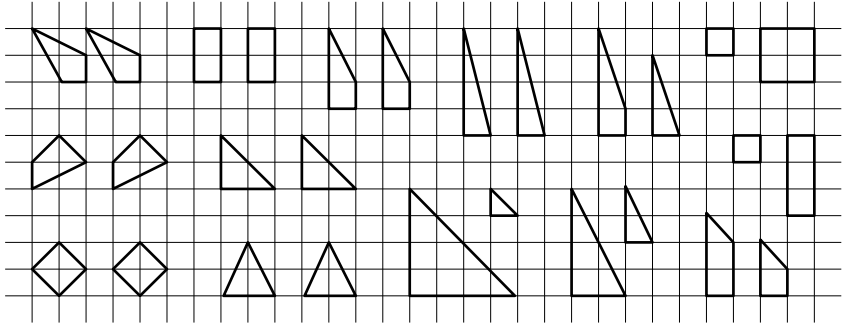}
\end{center}
\label{Mixed_area_4}
\caption{Maximal pairs of polygons of mixed volume $4$}
\end{figure}

\end{thm}

The ``only if'' part can be reduced to Theorem~\ref{mix1} by means of Example \ref{lucky2}, the classification is proved in Section \ref{Slattice}, and the ``if'' part follows from the observation that each of the listed 14 systems is explicitly solvable.

We cannot prove Conjecture 1 in full generality for $n>2$, nor can we classify tuples of lattice polytopes of mixed volume 4 in $\Z^n,\, n>2$. The rest of this section is devoted to solving these problems in certain special cases.

\begin{dfntn} A point $a\in A_1$ is said to be a {\it lucky} point of a tuple $A_1,A_2,\ldots,A_n\subset\Z^n$, if every tuple of faces $\Gamma_j$ of the convex hulls of $A_j$ such that $a\in \Gamma_1$, $\dim\Gamma_1>0$, and $\G= \sum_j\Gamma_j$ is a face of the convex hull of $\sum_j A_j$ satisfies the following two conditions:
\begin{enumerate}
\item if there exists $I\subsetneq\{1,2,\ldots,n\}$ such that $\dim\sum_{j\in I}\Gamma_j<|I|$, then $I$ can be chosen so that $1\notin I$.

\item if there is no such $I$, then the lattice distance from the affine hyperplane $L$ containing $\G$ to the set $\left(\sum_j A_j\right)\setminus L$
is smaller than the half of the mixed volume of the convex hulls of $A_1,A_2,\ldots,A_n$.
%the meaningful case is when a is vertex, but everything remains valid a forteriori if a is not!
\end{enumerate}
\end{dfntn}

\begin{exa}\label{lucky1} If the volume of the convex hull of an irreducible set $A\subset\Z^n$ is greater than 2, then every $a\in A$ is a lucky point of the tuple $(A,A\setminus\{a\},\ldots,A\setminus\{a\})$. Moreover, the mixed volume of the convex hulls of this tuple equals the volume of the convex hull of $A$.
\end{exa}

{\it Proof.} The equality of volumes easily follows from the Kouchnirenko---Bernstein formula.
For any collection $\Gamma_1,\Gamma_2,\ldots\Gamma_n$ in the definition of the lucky point, we have $\G_2=\G_3=\ldots=\G_n$ and $\G_1=\conv(\G_2\cup \{a\})$. Therefore, condition 1 is always satisfied: for $I\subsetneq\{1,2,\ldots,n\}$ such that $\dim\sum_{j\in I}\Gamma_j<|I|$, the same is true for $I'=I\setminus\{1\}\neq \emptyset$. Now we examine condition~2 in the case when $\G_1=F\ni a$ is a facet of $\conv A$: if $\Gamma_2$ is also a facet, then the lattice distance $d$ from $A\setminus F$ to the hyperplane containing $F$ is smaller than $\frac{1}{2} \Vol \conv A$. Note that we have $d\Vol F\leqslant\Vol\conv A$.
\begin{itemize}
\item Assume that $\Vol F=1$. Then $a$ is a vertex, and $\dim \sum_{j\in I} \G_j = \dim \G_2 = n-2 < |I|$, for $I=\{2,3,\ldots, n\}$.

\item Assume that $\Vol F=2$ with the equality $d\Vol F=\Vol\conv A$. This case is impossible: the latter equality implies that $A\setminus F$ consists of one point placed at lattice distance $d$ from the affine span of $F$. Since $\Vol A>2$, we have $d>1$, and thus $A$ is reducible.

\item Assume that $\Vol F>2$ or $\Vol F=2$ with the strict inequality $d\Vol F<\Vol\conv A$. In this case $d<1/2\Vol\conv A$. $\quad\Box$
 \end{itemize}

\begin{exa}\label{lucky2} For every reduced pair $A, B\subset\Z^2$ and every point $a\in A$, there exists $\tilde B\subset B$ such that $a$ is a lucky point of the pair $(A,\tilde B)$ and the mixed area of the convex hulls of $A$ and $\tilde B$ is the same as for $A$ and $B$.
\end{exa}

The proof is elementary and follows from the same considerations as the preceding one.

\begin{thm}\label{mix1}
Conjecture 1 is valid for every tuple with a lucky point. Moreover, assuming with no loss of generality that the lucky point is $0\in A_1$, the solution of a system of equations
\begin{equation}\label{System_whith_one_parameter}
\sum_{a\in A_1} c_{1,a}x^a=c,\, \sum_{a\in A_2} c_{2,a}x^a=\ldots=\sum_{a\in A_n} c_{n,a}x^a=0
\end{equation}
with generic coefficients $c_{j,a}\in\C$ cannot be expressed by generalized quadratures in terms of the right hand side $c$.
\end{thm}

In the preceding sentence, the word ``solution'' refers to the multivalued function $F_{c_{\AA}}\colon \C \to (\C\setminus 0)^n$ whose value at every point $c\in \C$ is the set of solutions of the system~(\ref{System_whith_one_parameter}).

\begin{thm}
\label{mix2}
Conjecture 1 is valid for $A_1,A_2\ldots,A_n$, if there exist two points of $A_1$ such that the segment connecting these points is not contained in the boundary of the convex hull of $A_1$.
\end{thm}

The proof is given in Section~\ref{Sproofs2} and is based on a certain generalization of the fact that the transitive group generated by transpositions is symmetric, which may be 
of independent interest.

\begin{rmrk} The roots of the general equation $\sum_j c_jx^j=0$ can be regarded as a multivalued function $R_d$ of the coefficients $c_0,\ldots,c_d$. A general system of equations is said to be {\it solvable by $N$-radicals}, if its solution can be expressed in terms of its coefficients using arithmetic operations and the functions $R_2,R_3,\ldots,R_N$. The proof of Theorem \ref{mix2} allows to strengthen its conclusion as follows: the general system of equations with $N$ roots, supported at a reduced tuple with an interior segment, is not solvable by $(N-1)$-radicals. This is because in Section~\ref{Sproofs2} we actually prove that the monodromy group of the system equals $S_N$, which implies non-solvability in $d$-radicals for $d<N$ (see \cite{topkhov}).

The other preceding theorems, in contrast, cannot be strengthened in this way, because our proof of non-solvability of the monodromy group in Section \ref{Smonodr} is not based on proving that the monodromy group is symmetric. In particular, we cannot prove that the general system of equations of degrees $d_1,\ldots,d_n$ is not solvable by $(d_1\cdot\ldots\cdot d_n-1)$--radicals (note that Theorem \ref{mix2} is not applicable to such systems for many tuples $d_1,\ldots,d_n$ such that $d_1\cdot\ldots\cdot d_n>4$).
\end{rmrk}

\begin{rmrk} Theorems~\ref{mix1} and \ref{mix2} are not enough to prove Conjecture 1 in dimension 3, i.e. there do exist reduced triples of lattice sets in $\Z^3$ with mixed volume greater than 4 and 
with no lucky points and interior segments.\end{rmrk}

\section{Classification of small lattice polytopes} \label{Slattice}

\hspace{3ex} {\it Proof of the classification presented in Theorem~\ref{Thm_equal_polytops}} is based on the notion of a circuit.

\begin{dfntn} A set $A\subset\Z^n$ is said to be {\it affinely dependent}, if there exist coefficients $c_a\in\R$ such that $\sum_{a\in A} c_a=1$ and $\sum_{a\in A} c_aa=0$. A {\it circuit} is a minimal (by inclusion) affinely dependent set.
\end{dfntn}
Every circuit $A$ can be uniquely decomposed into a disjoint union $A=A_+\sqcup A_-$ such that $A_+$ and $A_-$ are the sets of vertices of two simplices with a unique common interior point (this decomposition is induced by the signs of the coefficients $c_a$ in the unique affine dependence relation for $A$).

The lattice volume of the convex hull $\conv A$ equals the product of the lattice volumes of $\conv A_+$ and the projection of $\conv A_-$ along the affine span of $A_+$. Since the latter has an interior lattice point, its volume is not smaller than $|A_-|$. In particular, $\Vol\conv A\leqslant 4$ implies $|A_\pm|\leqslant 4$, so the dimension of the affine span of $A$, which equals $|A|-2$, is at most 6. This observation leads to the following classification:
\begin{lem} \label{classcirc} Every circuit $A$ such that $\Vol\conv A\leqslant 4$ coincides (up to an affine automorphism of the lattice) with one of the circuits from Theorem~\ref{Thm_equal_polytops}, or $\{(\pm 1,0), (0,\pm 1\}$, or $\{0,2,4\}$.
\end{lem}

Let now $A\subset\Z^n$ be an arbitrary irreducible set such that $\Vol\conv A\leqslant 4$. Among all circuits of maximal volume in $A$, choose a circuit of maximal cardinality $B\subset A$. This circuit is non-trivial, unless $A$ is the set of vertices of a simplex, and we can assume with no loss in generality that $0\in B$. We shall now classify all possible $A$ with a given circuit $B$, where $B$ runs over all circuits listed in Lemma \ref{classcirc}.

If $\Vol\conv B=4$, then $A$ is obviously an iterated standard cone over $B$ (up to an affine automorphism of $\Z^n$).

If $\Vol\conv B=3$, consider the image $A'$ of $A$ under the projection along the vector span of $B$. We need the following observation:
\begin{lem} \label{projl} For every irreducible set $M\subset\Z^m$ with $0\in M$, exactly one of the following possibilities takes place:

1) There exist disjoint simplices of total volume 3 or more with a common vertex 0, such that all of their vertices are in $M$.

2) $\conv M$ is an iterated standard cone over a segment of lattice length 2 or over a parallelogram of lattice area 2. In this case, there exist disjoint simplices of total volume 2 with a common vertex 0, such that all of their vertices are in $M$.

3) $M$ is the set of vertices of a unit simplex.
\end{lem}

Applying this lemma to $M=A'$, we conclude that only the third possibility can take place, otherwise the volume of $\conv A$ would be at least $3\cdot 2=6$. Since $A'$ is the set of vertices of a unit simplex, and $\Vol\conv A'=3<\Vol\conv A=4$, then at least one of the points in $A'\setminus \{0\}$ is the image of at least two points of $A$, whose difference $b$ is contained in the vector span of $B$. Thus, $A$ contains an iterated standard cone over the set $B'=B\times\{0\}\cup\{0,b\}\times\{1\}$. Since $\Vol\conv B'\geqslant 4$, we conclude that $A$ is an iterated standard cone over $B'$, and $\Vol\conv B'=4$. For each of $B$ from Lemma \ref{classcirc}, it is an elementary-geometric problem to classify all suitable $b$ such that $\Vol\conv B'=4$.

Finally, if $\Vol\conv B=2$, we also consider the image $A'$ of $A$ under the projection along the vector span of $B$. If $A'$ is the set of vertices of a unit simplex, then, similarly to the case $\Vol\conv B=3$, at most two points of $A'\setminus\{0\}$ are images of more than one point of $A$, and $A$ is an iterated standard cone over $B'=B\times\{(0,0)\}\cup B_1\times\{(0,1)\}\cup B_2\times\{(1,0)\}$, where $B$ is either $\{0,1,2\}$ or $\{0,1\}\times\{0,1\}$, and the sets $B_1, B_2\subset\Z^2$ are easy to classify for $\Vol\conv B'\leqslant 4$. Finally, if $A'$ is not the set of vertices of a unit simplex, then, by Lemma \ref{projl}, we observe that $\Vol\conv A'=2$, and every point of $A'\setminus\{0\}$ is the image of a unique point of $A$, so $A$ can be reconstructed uniquely from $B$ and $A'$. $\quad\Box$

{\it Proof of the classification presented in Theorem~\ref{Thm_dimension_two}} is based on the following well known formula for the mixed area of polygons:

\begin{lem} Let $L_i(x,y)<c_i$ be the minimal set of inequalities describing a lattice polygon $A$, where $L_i:\Z^2\to\Z$ are surjective linear functions. Denote the lattice length of the edge of $A$ on the line $L_i(x,y)=c_i$ by $a_i$, and the maximal value of $L_i$ on a polygon $B$ by $b_i$, then the mixed volume of $A$ and $B$ equals $\sum_i a_ib_i$.
\end{lem}

By the Aleksandrov--Fenchel inequality $\MV(A,B)^2\geq\Vol(A)\Vol(B)$, assuming $\Vol(A)\leqslant\Vol(B)$, we have $\Vol(A)\leqslant 4$. The classification of all $A$ of area $\leq 4$ is given above. Let $A$ be one of these polygons, and define the functions $L_i$ and the numbers $a_i$ as in the lemma above. Then there are finitely many tuples of non-negative integers $(b_i)$ such that  $\sum_i a_ib_i\leqslant 4$. Up to a parallel translation, every $B$ such that $\MV(A,B)\leqslant 4$ is contained in the intersection of the half-planes $L_i<b_i$ for one of such tuples $(b_i)$. $\quad\Box$

\section{Monodromy of general systems of equations: preliminaries}\label{Smonodr0}

The proof of our results is based on the search for general systems of equations whose monodromy group is solvable. This condition is equivalent to the solvability of the system itself by the following topological version of Galois theory (see e.g. \cite{topkhov}):

\begin{thm} \label{topgalois0} Let $C\subset\C^2$ be an algebraic curve, $(x,y)$ the standard coordinates, and $\pi$ the restriction of $x$ to $C$. The multivalued function $y\circ\pi^{-1}$ can be expressed in generalized quadratures if and only if the monodromy group of the branched covering $\pi$ is solvable.
\end{thm}

We shall need the following version of this fact.

\begin{prop} \label{topgalois1} Let $C\subset\C^n$ be a curve and $f:\C^n\to\C$ a polynomial. The coordinates of the points $x\in C,\, f(x)=c$, can be expressed in generalized quadratures in terms of $c$ if and only if the monodromy group of the branched covering $f:C\to\C$ is solvable.
\end{prop}
{\it Proof.} Let $C'\subset\C^2$ be the image of the map $(f,x_j):C\to\C^2$, where $(x_1,\ldots,x_n)$ are the standard coordinates in $\C^n$. Denoting the restriction of $f$ to $C'$ by $\pi'$, we observe that the monodromy group of $\pi'$ is the same as for $f:C\to\C$, so Theorem \ref{topgalois0} for $C'$ and $\pi'$ gives the statement of the proposition for the $j$-th coordinate. $\quad\Box$

We thus start with counting the monodromy for a general system of equations. For $i=0,\ldots,m$, let $A_i$ be a finite set in $\Z_{\geqslant 0}\times\Z^m,\, 0\in A_i$, and let $\C^{A_i}$ be the space of Laurent polynomials of the form $f(t,x)=\sum_{a\in A_i} c_a t^{a_0} x_1^{a_1}x_2^{a_2}\ldots x_m^{a_m},\, t\in\C,\, x\in\CC^m$. Denote the projection $\Z_{\geqslant 0}\times\Z^m\to\Z_{\geqslant 0}$ by $\mu$ and the (non-empty) intersection $A_i\cap\ker\mu$ by $B_i$. 
Then, for a generic tuple $(f_0,f_1,\ldots,f_n)\in\prod_{i=0}^m\C^{A_i}$, we choose $c_0$ to be the value $f_0(0,x)$ for one of finitely many $x$ such that $f_1(0,x)=f_2(0,x)=\ldots=f_n(0,x)=0$, and will be interested in the monodromy of the finite set $\{f_0=c,\, f_1=f_2=\ldots=f_n=0\}$ as $c$ runs a small circle around $c_0$. The answer can be formulated in terms of the numbers $m_i=\min\mu|_{A_i\setminus B_i}$ and can be extracted from the following important fact.

\begin{lem} \label{keylemma1} In the setting of the preceding paragraph, assume that the dimension of the vector span of $\bigcup_{i\in I} B_i$ is at least $|I|$ for every $I\subsetneq\{0,1,\ldots,m\}$. Then, choosing a generic $(f_0,f_1,\ldots,f_m)\in\prod_{i=0}^m\C^{A_i}$ among all tuples such that $f_0(0,x)=f_1(0,x)=\ldots=f_m(0,x)=0$ for a given $x\in\CC^n$, the intersection multiplicity of the surfaces $\{f_i=0\},\, i=0,\ldots,m$, at the point $(0,x)$ is well defined and equals the minimum of $m_i,\, i=0,\ldots,m$.
\end{lem}

The proof is based on two following obvious observations.

\begin{lem} \label{contact1} If smooth hypersurfaces $H_1,\ldots,H_m$ in $\C^{m+1}$ are mutually transversal, and a smooth curve, passing through a point $x\in\bigcap_i H_i$, has an order $m_i$ contact with $H_i$ at $x$, then it has an order $\min_i m_i$ contact with the curve $\bigcap_i H_i$ at $x$.
\end{lem}

\begin{lem} \label{contact2} Let $L$ be a line parallel to the $x_0$ axis of the space $\C^{m+1}$ with the standard coordinates $(x_0,x_1,\ldots,x_m)$, and let $y$ be the intersection $L\cap \ker x_0$. Assume that a smooth function $f:\C^{m+1}\to\C$ has no critical point at $y$, but its restriction to $L$ has a root of order $m$ at $y$. Assume that a smooth curve $C\subset\C^{m+1}$ has a contact of order $\geqslant m$ with $L$ at $y\in\ker x_0$. Then, for generic $\alpha\in\C$, the restriction of $f(\alpha x_0,\, x_1,\, x_2,\ldots,x_{m})$ to $C$ has a root of order $m$ at $y$.
\end{lem}

{\it Proof of Lemma \ref{keylemma1}.} With no loss in generality, we assume that $m_0\leqslant m_i$ for $i>0$. Since the dimension of the vector span of $\bigcup_{i\in I} B_i$ is at least $|I|$ for every $I\subsetneq\{0,1,\ldots,m\}$, the hypersurfaces $f_i=0,\, i=1,2,\ldots,m$ and $t=0$ are transversal at $(0,x)$ by genericity (see e.g. \cite{sturmf94}). Also by genericity, the order of contact of $f_i=0$ and the line $L=\{(t,x)\,|\,t\in\R\}$ equals $m_i$. Thus, by Lemma \ref{contact1}, the curve $f_1=f_2=\ldots=f_m=0$ has a contact of order $m'=\min_{i=1}^m m_i$ with $L$. Then the restriction of $f_0$ to this curve has a root of order $m_0$, otherwise $f_0$ is not generic by Lemma \ref{contact2}. The sought intersection number equals the order of this root. $\quad\Box$

In order to drop the linear-algebraic assumption in Lemma \ref{keylemma1}, we need the following notation. To finite sets $B_i\ni 0,\, i=0,1,\ldots,m$, in $\Z^m$, assign a number $d_B$, a set $I_B\subset\{0,1,\ldots,m\}$ and a rational subspace $L_B\subset\R^m$ as follows. For every $I\subset\{0,1,\ldots,m\}$, define $L_I$ as the vector span of the union $\bigcup_{i\in I} B_i$. If there exists $I\subset\{1,2,\ldots,m\}$ such that $\dim L_I<|I|$, then we set $d_B=0$, and $I_B$ and $L_B$ are not defined. Otherwise, among all $I\subset\{0,1,\ldots,m\}$ such that $\dim L_I<|I|$, there exists the minimal one by inclusion (see \cite{sturmf94}). We denote this minimal set by $I_B$, define $L_B$ as $L_{I_B}$, and define $d_B$ as the product $($the index in $L_B$ of the sublattice generated by $\bigcup_{i\in I_B} B_i)\cdot($the mixed volume of the convex hulls of the images of $B_i,\, i\notin I_B$, under the projection $\Z^m\to\Z^m/L_B)$.

\begin{lem}[\cite{sombra15}, \cite{khov15}]\label{lemmad} For generic polynomials $g_i\in\C^{B_1}$ and a number $c_0$ such that the system of equations $g_0(x)=c_0,\, g_1(x)=\ldots=g_m(x)=0$ is consistent, this system has $d_B$ solutions.
\end{lem}

\begin{prop} \label{localmonodr}
In the setting of the paragraph preceding Lemma \ref{keylemma1}, for a generic tuple \linebreak $(f_0,f_1,\ldots,f_m)\in\prod_{i=0}^m\C^{A_i}$, the permutation of the finite set $\{f_0=c,\, f_1=f_2=\ldots=f_n=0\}$ as $c$ runs a small circle around $c_0$ consists of $d_B$ cycles of length $\min_{i\in I_B} m_i$ each.
\end{prop}
{\it Proof.} The restriction of $f_0$ to the curve $f_1=f_2=\ldots=f_n=0$ has finitely many simple roots outside the hyperplane $t=0$, and $d_B$ roots of order $\min_{i\in I_B} m_i$ each at this hyperplane. The number of roots is computed by Lemma \ref{lemmad}. The computation of the order of each of these roots can be obviously reduced to the case $L_B=\R^m$ and then done by Lemma \ref{keylemma1}. $\quad\Box$.

\begin{exa} For $f(x,y)=x^a+\alpha y^b-1$ and $g(x,y)=x^c+y^d-1$ with generic $\alpha$, the permutation of the roots of the system $f=\varepsilon,\, g=0$, as $\varepsilon$ runs around 0, consists of $\gcd(a,c)$ cycles, each of length $\min(b,d)$.
\end{exa}

\section{Monodromy of general systems of equations: the answer}\label{Smonodr}

In the preceding section, we have explained how the solvability of a system is related to the solvability of its monodromy, and described the monodromy of a general system at a single branching point. In this section, we completely describe the monodromy of a general square system (c.f. \cite{sturmf11}), and also outline another way to obtain this answer. This alternative way is less elementary, but more powerful (in particular, it extends to non-square systems).

For a reduced tuple $A=(A_1,A_2,\ldots,A_n),\, A_j\subset\Z^n$, let $\mathcal{B}\subset \C^{\AA}=\C^{A_1}\oplus\C^{A_2}\oplus \ldots\oplus\C^{A_n}$ be the {\it bifurcation set}, i.e. the closure of all $c\in\C^{\AA}$ such that the system
$$
\sum_{a\in A_j} c_{j,a} x^a=0 \mbox{ for } j=1,2,\ldots, n\eqno{(*)}
$$
has less solutions than we expect for systems defined by generic $c\in\C^{\AA}$ (i.e. less than the mixed volume of the convex hulls of $A_1,A_2,\ldots,A_n$).

The irreducible components of $\mathcal{B}$ are in one to one correspondence with the {\it essential facings} of the tuple $A_1,A_2,\ldots,A_n$ (defined below), and we shall describe the cycle type of the permutation of the solutions of the system $(*)$ as $c$ travels around each of these components.

\begin{dfntn}
The {\it codimension} of a tuple of sets $B_1,B_2,\ldots,B_k\subset\Z^n$ is the difference $k-\dim(\conv \sum_j B_j)$.
A tuple of subsets $F_j\subset A_j,\, j=1,2,\ldots,n,$ is said to be a {\it face} of $A_1,A_2,\ldots,A_n$, if the following equivalent conditions are satisfied:

1) $F_j=\Gamma_j\cap A_j$, where $\Gamma_j$ is a face
of the convex hull of $A_j$, and the Minkowski sum $\sum_{j}\Gamma_j$ is a face of the convex hull of $\sum_{j} A_j$.

2) There exists a linear function $L:\Z^n\to\Z$, whose restriction to $A_j$ attains its maximum at $F_j$ (i.e. $L(x)<L(y)=\max L(A_j)$ for $y\in F_j$ and $x\in A_j\setminus F_j$).

A subtuple $S$ of a face of $A_1,A_2,\ldots,A_n$ is said to be an {\it essential facing} of $A_1,A_2,\ldots,A_n$, if its codimension is $1$, and the codimension of every proper subtuple of $S$ is at most 0.
\end{dfntn}

For an essential facing $S=(S_{l_1},S_{l_2},\ldots,S_{l_p})$ of the tuple $A=(A_1,A_2,\ldots,A_n)$, $S_{l_j}\subset A_{l_j}$,
let $R_S$ be the closure of all $c\in\C^{\AA}$ such that the equations $\sum_{a\in S_{l_j}} c_{l_j,a} x^a=0 \mbox{ for } j=1,2,\ldots,p$ are compatible. Let $D_{\AA}$ be the set of all $c$ such that the equations $(*)$ have a multiple solution.
\begin{prop}[\cite{e13}] \label{monodrn1} If the tuple $A$ is reduced, then the sets $D_A$ and $R_S$ for all of the essential facings $S$ are pairwise different irreducible hypersurfaces, and their union $D_{\AA}\cup \bigcup_S R_S$ equals the bifurcation set $\mathcal{B}$.
\end{prop}

We now describe the cycle type $T_S$ of the permutation of the solutions of $(*)$ as $c$ runs around $R_S$. We shall encode the type of a permutation with $a_i$ cycles of length $i$ for $i\in \N$ as $\sum_i a_ie_i\in\Z^\N$, where $e_1,e_2,\ldots$ is the standard basis in $\Z^\N$.

Let $S=(S_{l_1},S_{l_2},\ldots,S_{l_p})$ be an essential facing of the tuple $A=(A_1,A_2\ldots,A_n)$ contained in its codimension 1 face $F=(S_1,S_2,\ldots,S_n)$. Permuting and shifting $A_1,A_2,\ldots,A_n\subset\Z^n$, we can provide that $S=(S_1,S_2,\ldots,S_p)$, and $S_j\ni 0$ for $j=1,2,\ldots,n$. Under these assumptions, let $L(S)$ be the lattice generated by $S_1\cup S_2\cup\ldots\cup S_p$, let $\bar L(S)$ be its saturation $L(S)\otimes{\mathbb Q}\cap\Z^n$, denote the index $|\bar L(S)/L(S)|$ by $i(S)$, and the projection $\Z^n\to\Z^n/\bar L(S)$ by $\pi_S$. There exists a unique primitive covector $l:\Z^n\to\Z$ whose restriction to $A_j$ attains its minimum on $S_j$ for $j=1,2,\ldots,n$ (by our assumptions, this minimum equals 0). We denote the minimum of $l$ on $\bigcup_{j=1}^p (A_j\setminus S_j)$ by $l(F)\in\N$, and the mixed volume of the convex hulls $\pi_S(S_{p+1}),\ldots,\pi_S(S_{n})$ in $\ker l/\bar L(S)$ by $v(F)$.
\begin{prop} \label{monodrn2} 1) If the tuple $A$ is reduced, then, as $c$ runs around
$D_A$, two solutions of the system $(*)$ permute.

2) For every essential facing $S$, the corresponding permutation type $T_S$ of a loop around $R_S$ equals $i(S)\sum_F v(F) e_{l(F)}$, where $F$ runs over all codimension 1 faces of $A$ such that $S$ is a subtuple of $F$.
\end{prop}
{\it Proof.} Part 1 is obvious, and Part 2 rephrases Proposition \ref{localmonodr}. $\quad\Box$.

We shall also obtain Propositions~\ref{monodrn1} and~\ref{monodrn2} in a less straightforward manner as a special case of the computation of monodromy $\zeta$-functions for a system of $k$ equations in $n$ variables with arbitrary $k\leqslant n$ (Theorem~\ref{zeta1} below). Although the case $k<n$ is not necessary for the purpose of this paper, it may be of independent interest, so we switch to the general case $k\leqslant n$ till the end of this section.

Consider the graph of the solution for the general system of equations with the Newton polytopes $A_1,A_2,\ldots,A_k$ in $\Z^n$:
$$
\Gamma=\{ (x,c) \, |\, \sum_{a\in A_j} c_{j,a} x^a=0 \mbox{ for } j=1,2,\ldots,k \} \subset \CC^n\times\C^{\AA}.
$$

Let $\mathcal{B}$ be the bifurcation set of the solution, i.e. the set of all points in $\C^{\AA}$, at which the projection of $\Gamma$ to $\C^{\AA}$ fails to be a locally trivial fibration. Recall the description of the irreducible components of $\mathcal{B}$ obtained in \cite{e13}.

\begin{dfntn}
A subtuple of $B_1,B_2,\ldots,B_k\subset\Z^n$
is said to be {\it important}, if it cannot be extended to a subtuple of higher codimension.
\end{dfntn}

\begin{dfntn}
For any subsets $(B_{i_1}\subset A_{i_1},\ldots,B_{i_p}\subset A_{i_p})$, define $i_{B_{i_1},\ldots,B_{i_p}}$ as the index of the sublattice of $\Z^{n+1}$, generated by the sets $B_{i_j}\times\{1\},\, j=1,\ldots,p$. Define the algebraic set $D_{B_{i_1}*\ldots*B_{i_p}}$ as the closure of all $c\in\C^{A_1}\oplus\ldots\oplus\C^{A_k}$ such that 0 is a critical value of the polynomial $\sum_j \lambda_j \sum_{a\in B_{i_j}} c_{i_j,a} x^a$ of the variables $\lambda\in\CC^p$ and $x\in\CC^n$.
\end{dfntn}

Note that the algebraic set $D_{B_{i_1}*\ldots*B_{i_p}}$ is irreducible.

\begin{prop} \label{pbif} Assume that a tuple $(A_1,\ldots,A_k)$ does not contain subtuples of positive codimension.

1) The irreducible components of the bifurcation set $\mathcal{B}$ are the hypersurfaces of the form $D_{B_{i_1}*\ldots*B_{i_p}}$, where $(B_{i_1}\subset A_{i_1},\ldots,B_{i_p}\subset A_{i_p})$ runs over all important subtuples of faces of $(A_1,\ldots,A_k)$.

2) At a generic point $c\in D_{B_{i_1}*\ldots*B_{i_p}}$, the Euler characteristic of the fiber of the projection  $\Gamma\to\C^{A_1}\oplus\ldots\oplus\C^{A_k}$ drops by the number $i_{B_{i_1},\ldots,B_{i_p}}c_{A_1,\ldots,A_k}^{B_{i_1},\ldots,B_{i_p}}$ as compared to a generic fiber of this projection.
\end{prop}
See the remark after Theorem~\ref{zeta1} below or Definition 2.27 in \cite{e13} for the definition of the number $c_{A_1,\ldots,A_k}^{B_{i_1},\ldots,B_{i_p}}$, and see Proposition 1.11 and Corollary 2.29 in \cite{e13} for the proof.

A loop around an irreducible component of the bifurcation set $\mathcal{B}$ induces automorphisms $h_i:H^i\to H^i$ in the cohomology of the generic fiber of the projection $\Gamma\to\C^{A_1}\oplus\ldots\oplus\C^{A_k}$. We shall describe the corresponding $\zeta$-function, i.e. the rational function $\zeta(t)=\prod_i \det(h_i - t\cdot\Id)^{(-1)^i}$.

For a subtuple $B=(B_{i_1}\subset A_{i_1},\ldots,B_{i_p}\subset A_{i_p})$ of a proper face of $(A_1,\ldots,A_k)$, denote the projection of $\Z^n$ along the affine span of $B_{i_1}+\ldots+B_{i_p}$ by $\pi_B:\Z^n\to L_B$, the sum $A_{i_1}+\ldots+A_{i_p}$ by $B'$, and the convex hull of $(A_{i_1}\cap\Z^n\setminus B_{i_1})+\ldots+(A_{i_p}\cap\Z^n\setminus B_{i_p})$ by $B''$. The tuple of the images $\pi_B A_1,\ldots, \pi_B A_k, \pi_B B', \pi_B B''$ has finitely many faces $C=(C_1,\ldots,C_k, C', C'')$ such that $C'\in B'$ is the vertex $\pi_B(B_{i_1}+\ldots+B_{i_p})$, and the sum $C_1+\ldots+C_k+C''$ has codimension 1 in $L_B$. We denote the set of all such faces by $F_B$, and for every such face $C\in F_B$ define two integer numbers:

\vspace{1ex}

\noindent $h_C=|l(C'')-l(C')|$, where $l$ is the surjective linear function $l:L_B\to\Z$ whose restriction to $C_1+\ldots+C_k+C''$ is constant.

\vspace{1ex}

\noindent $m_C=\sum\MV\left(\underbrace{C'',\ldots,C''}_{m_0-1},\underbrace{C_1,\ldots,C_1}_{m_1},\ldots,\underbrace{C_k,\ldots,C_k}_{m_k} \right)$, where the sum is taken over all positive integers $m_0,\ldots,m_k$ that sum up to $\dim L_B$.

\begin{thm}\label{zeta1}
1) The $\zeta$-function of monodromy of the projection $\Gamma\to\C^{A_1}\oplus\ldots\oplus\C^{A_k}$, restricted to a loop around $D_{A_1*\ldots*A_k}$, equals $$(1-t^2)^{i_{A_1,\ldots,A_k}}.$$

2) If $B=(B_{i_1}\subset A_{i_1},\ldots,B_{i_p}\subset A_{i_p})$ is an important subtuple of a proper face of $(A_1,\ldots,A_k)$, then the $\zeta$-function of monodromy of the projection $\Gamma\to\C^{A_1}\oplus\ldots\oplus\C^{A_k}$, restricted to a loop around $D_{B_{i_1}*\ldots*B_{i_p}}$, equals $$\prod_{C\in F_B} (1-t^{h_C})^{i_Bm_C}.$$
\end{thm}

Note that the degree of the latter rational function equals $i_{B_{i_1},\ldots,B_{i_p}} c_{A_1,\ldots,A_k}^{B_{i_1},\ldots,B_{i_p}}$, i.e. Theorem~\ref{zeta1} implies Part 2 of Proposition \ref{pbif} that cites Theorem~2.36 in \cite{e13}. The proof is also the same as for Theorem~2.36 in \cite{e13}, extending the computation of Milnor numbers of the fibers by methods of \cite{e10} to the computation of monodromy $\zeta$-functions of the fibers by methods of \cite{mt2}. The rest of the proof of Theorem~2.36 also literally extends to monodromy $\zeta$-functions, because they enjoy the same additivity properties as the Euler characteristic (see e.g. \cite{smg}).

\section{Variations on Ritt's lemma and the proof of Theorem~\ref{mix1}} \label{Smix1}

\begin{thm}[\cite{jones}]\label{jordan}
If a primitive subgroup of $S_n$ contains a cycle of length at most $n-3$, then it equals $A_n$ or $S_n$.
\end{thm}

We shall use this extension of the classical Jordan theorem to prove Theorem~\ref{mix1}. Choose generic $c_{i,a}\in\C,\, a\in A_i$, set $c_{0,1}=0$ and denote $\sum_{a\in A_i} c_{i,a}x^a$ by $f_i(x)$. Also denote the curve $f_2=\ldots=f_n=0$ by $X\subset\CC^n$, and the mixed volume of the convex hulls of $A_1,\ldots,A_n$ by $d$.
We assume that $d>4$ and wish to prove that the monodromy group $M$ of the degree $d$ branched covering $f_1:X\to\C$ is not solvable; then the solution of the system $f_1(x)=c,\, f_2(x)=\ldots=f_n(x)=0$ cannot be expressed in terms of $c$ by generalized quadratures by Proposition \ref{topgalois1}.

Note that $M$ is transitive, because the reducedness of the tuple $(A_1,\ldots,A_n)$ implies that $X$ is connected by the following lemma:

\begin{lem}[\cite{khov15}] \label{lconnect} Let $m$ be the maximum of the codimensions of all the subtuples in a tuple $B_1,\ldots,B_k\subset\Z^d$, and consider a generic complete intersection $Z$ given by the equations $\sum_{a\in B_i} c_{i,a}z^b=0,\, i=1,\ldots,k$.

1) If $m$ is positive, then $Z$ is empty.

2) If $m$ equals 0, then one can choose the maximal (by inclusion) subtuple $B_{i_1},\ldots,B_{i_p}$ of codimension 0, and the number of connected components of $Z$ equals the $p$-dimensional mixed volume of the convex hulls of $B_{i_1},\ldots,B_{i_p}$.

3) If $m$ is negative, then $Z$ is connected.
\end{lem}

The notion of a lucky point and the description of monodromy in the preceding section imply that the bifurcation set of $f_1:X\to\C$ consists of the point 0, whose monodromy is arbitrarily complicated, and finitely many other points, whose monodromy is each a single cycle of length smaller than $d/2$.

Since $M$ contains transpositions, corresponding to the critical points of $f_1$, then, by the Jordan theorem, the primitivity of $M$ implies $M=S_d$ or $A_d$, i.e. $M$ is not solvable. Thus, it remains to study the case of imprimitive $M$. In this case, the branched covering $f_1:X\to\C$ splits into a non-trivial composition $X\stackrel{g}{\to} Y\stackrel{h}{\to} \C$, where the maps $g$ and $h$ induce the structure of a Riemann surface on $Y$. This was first noticed by Ritt in \cite{ritt} under the assumption that the genus of $X$ is 0, but remains valid for arbitrary genus. Moreover, we can assume with no loss in generality that the monodromy group $M'$ of the branched covering $h:Y\to\C$ is primitive, otherwise we could decompose $h$ in the same way. Furthermore, $M'$ is a quotient of $M$, so non-solvability of $M'$ implies the same for $M$. Thus, it remains to consider the case of solvable $M'$.

Since $M'$ is primitive and solvable, the monodromy of a critical value of $h$ cannot be a cycle of length smaller than the half of the degree $d'$ of $h$: otherwise  the degree of $h$ is greater than 4, the length of the cycle is at least by 3 smaller, and $M'$ is not solvable by Theorem~\ref{jordan}.

On the other hand, the monodromy of a non-zero critical value of $h$ cannot be a cycle of length at least $d'/2$ or not a cycle, otherwise the monodromy of the same critical value of $f_1$ would be a cycle of length at least $d/2$ or not a cycle.

We conclude that $h:Y\to\C$ has no critical values besides 0, i.e. $Y=\C$, and $h(z)=z^{d'}$, i.e. $f=g^d$. This implies that the generic complete intersection curve $f_1-y^d=f_2=\ldots=f_n=0$ in $\CC^{n+1}$ has $d$ connected components, which contradicts Part 3 of Lemma \ref{lconnect}.

\section{Generating $S_n$ by arbitrary disjoint permutations} \label{Sproofs2} 

In this section, we prove a fundamental lemma stating that, under some appropriate conditions on a subset of permutations, the permutation group is generated by this subset. The proof of Theorem~\ref{mix2} is essentially based on this result.

\begin{dfntn}
We call a set of permutations $a_1,a_2,\ldots, a_t\in S_n$ {\it disjoint}, if for each $x\in \{1,2,\ldots, n\}$ there is at most one $j$ such that $a_j(x)\neq x$.
\end{dfntn}

\begin{lem}\label{Permutations}
Let a transitive subgroup $G\subset S_n$ be generated by a subset $\Sigma\subset S_n$ that contains disjoint non-trivial permutations $a_1,a_2,\ldots, a_t$ and whose other elements are transpositions. Assume $a_1$ has a fixed point $x$, i.e. $a_1(x)=x$. Then $G=S_n$.
\end{lem}

\begin{rmrk}
In the case $t>1$, the existence of fixed point $x$ is always satisfied.
\end{rmrk}

\begin{prf-2}
Each permutation $a_j$ is the product of its non-trivial cycles $c_{j,1}, c_{j,2}, \ldots, c_{j,k_j}$, i.e. $a_j= c_{j,1} c_{j,2} \ldots c_{j,k_j}$. The set of elements a cycle $c_{j,k}$ permutes is called a {\it carousel}. If an element $x\in \{1,2,\ldots, n\}$ does not belong to any carousel, we assume that it forms a virtual carousel consisting of the only point $x$. Since permutations $a_j$ (called {\it attractions}) are disjoint, all carousels are pairwise non-intersecting and form equivalence classes of elements. Consider the non-directed graph $\G=(V,E)$ whose vertices are carousels and two carousels $v_1,v_2\in V$ are connected, if there exists $x_j\in v_j$ such that the transposition $(x_1,x_2)$ belongs to $\Sigma$. Since $G$ is transitive, $\Gamma$ is connected. Consider an arbitrary spanning tree $T=(V,E')$ of the graph $\Gamma$. For each $w=(v_1, v_2)\in E'$, choose a transposition $a_w=(x_1, x_2)\in \Sigma$ such that $x_j\in v_j$. We further prove that the group $G'$ generated by $\Sigma'=\{a_j\mid j=1,2,\ldots, t\}\cup \{a_w\mid w\in E'\}$ coincides with $S_n$. It would be sufficient, since $G'\subset G$. For convenience, we consider any permutation $a\in S_n$ and are aiming to find $b\in G'$ such that $ba=e$ is the trivial permutation.

There exist $w=(v_1,v_2)\in E'$ such that $v_1$ is fixed under $a_1$ and $v_2$ is not. Informally, we consider $w$ as the root edge of the tree $T$. We prove, by induction on $k$, that there exist a subtree $T_k=(V_k,E'_k)$ of the tree $T$ and a permutation $b_k\in G'$ such that $v_1,v_2\in V_k$, $|V_k|\leq |V|-k$, and $b_k a(x)=x$ for all elements of the carousels $v\in V\setminus V_k$. Assume this to be proved for some integer $k<n-2$. Consider an arbitrary leaf $v_0\in V_k\setminus\{v_1,v_2\}$ of the tree $T_k$, $v_0=\{x_1,x_2,\ldots, x_l\}$. Take an arbitrary $x=x_j$. There exists a path $u_1, u_2,\ldots, u_{n_1}=v_1$ from the carousel $u_1\ni y:=b_ka(x)$ to the carousel $v_1$ and a path $u_{n_1+1}, u_{n_1+2},\ldots, u_{n_1+n_2}$ from $v_2=u_{n_1+1}$ to $v_0=u_{n_1+n_2}$ in the tree $T_k$. There exist integer $d_1,d_2,\ldots, d_{n_1+n_2}$ such that $b''_{x}(y)=x$, where $b''_x=a_{u_{n_1+n_2}}^{d_{n_1+n_2}} a_{n_1+n_2-1,n_1+n_2}\ldots a_{u_3}^{d_3} a_{23} a_{u_2}^{d_2} a_{12} a_{u_1}^{d_1}$, $a_{u_j}\in \Sigma'$ denotes the attraction that moves the carousel $u_j$ and $a_{j,j+1}$ is a transposition of two elements of carousels $u_j, u_{j+1}$. We aim to build $b_{k+1}$ in the form $b_{k+1} = \ldots b''_{x}b_k$, however, $b''_{x}$ may move some elements of $\bigcup_{v\in V\setminus V'} v$, which is not appropriate. Denote by $I=\{j\mid a_{u_j}=a_1\}$ the subset of carousels whose attraction is $a_1$. Consider
\begin{multline*}
b'_{x}:= a_{u_{n_1+n_2}}^{d_{n_1+n_2}} a_{n_1+n_2-1,n_1+n_2}\ldots a_{n_1-1, n_1} \left (\prod_{j\in I} a_{u_{j}}^{-d_{j}}\right ) a_{u_{n_1}}^{d_{n_1}} a_{n_1, n_1+1}\\ \left (\prod_{j\in \{1,2,\ldots, n_1+n_2\}\setminus I} a_{u_{j}}^{-d_{j}}\right ) a_{u_{n_1+1}}^{d_{n_1+1}} a_{n_1+1, n_1+2} \ldots a_{23} a_{u_2}^{d_2} a_{12} a_{u_1}^{d_1}.
\end{multline*}

The carousel $u_{n_1}$ is fixed under any attraction $a_{u_j}$ with $j\in I$, and the carousel $u_{n_1}$ is fixed under any attraction $a_{u_j}$ with $j\in \{1,2,\ldots, n_1+n_2\}\setminus I$, since it can be moved only by $a_1$. Therefore, $b'_{x}(y)=b''_{x}(y)=x$. On the other hand, $b'_{x}(z)=z$ for any $z\in \bigcup_{v\in V\setminus V_k} v$, since any such $z$ is fixed under $a_{j, j+1}$, $j=1,2,\ldots, n_1+n_2-1$. It is easy to see that $V_{k+1}:=V_{k}\setminus \{v_0\}$ and $b_{k+1}:=b'_{x_1}b'_{x_2}\ldots b'_{x_l}b_k$ satisfy the required conditions.

For $k=n-2$, only two carousels $v_1,v_2$ are not put in order by $b_{k}$. It is easy to see that the transposition of any two elements $x_1\in v_1$, $x_2\in v_2$ can be generated by $a_{v_1}$, $a_{v_2}$, and $a_{(v_1,v_2)}$. Therefore, $(b_k a)^{-1}$ belongs to $G'$ as well as any permutation supported at $v_1\cup v_2$.
\end{prf-2}

\section{Proof of Theorem~\ref{mix2}} \label{Smix2}

Consider a generic system $\ff=(f_1,f_2,\ldots, f_n)\in \C^{\AA}:= \C^{A_1}\oplus \C^{A_2}\oplus \ldots \oplus \C^{A_n}$. Consider the branched covering $\p_{\kk}\colon X = \{f_2=f_3=\ldots=f_n=0\} \to \C$ given by $\pi(\xx)=f_1(\xx)/\xx^{\kk}$. For a generic set $\{f_j\}$, a generic fiber $F=\p_{\kk}^{-1}(c)$ consists of $D$ points, where $D=\prod A_j$. Consider the group of permutations $S_{D}$ on $F$ and the monodromy group $G_{\kk}\subset S_{D}$ of the covering $\p_{\kk}$.

Since $S_D$ is not solvable for $|D|>4$, Proposition~\ref{topgalois1} implies Theorem~\ref{mix2} due to the following:

\begin{lem}\label{lem-main}
Assume there exist two integer points $\kk_0,\kk_1\in A_1$ such that segment $[\kk_0,\kk_1]$ is not contained in the boundary of $A_1$. If $A$ is reduced, then at least one of  $G_{\kk_0}, G_{\kk_1}$ coincides with~$S_D$.
\end{lem}

The rest of the paper is devoted to the proof of Lemma~\ref{lem-main}. As before, we use the following corollary of Lemma \ref{lconnect}:

\begin{lem}\label{transitivity}
If $\AA$ is reduced, then $X$ is connected.
\end{lem}

In what follows, we use the following notation for convenience. For a system $\AA=(A_1, A_2, \ldots, A_k)$ of subsets $A_j\in \Z^n$ whose Minkowski sum spans a $k$-dimensional sublattice $L\subset \Z^n$, we use $\prod \AA = \prod_{j=1}^k A_j$ for the mixed volume of the convex hulls $\conv A_1, \conv A_2,\ldots, \conv A_k$ in terms of the volume form in $L$. The set of primitive covectors $\a\colon \Z^n\to \Z$ is denoted by $\ZZ$.  For a set $A\subset \Z^n$ and a covector $\a\in \ZZ$, we use $A^\a$ for the face of~$A$, where $\a|_A$ attains its maximal value: $A^\a:= \{x\in A\mid \a(x)=\max \a|_A\}$. Moreover, we use $l_A(\a)$ for $\max \a|_A$. The proof of Lemma~\ref{lem-main} utilizes the following well-known equality for the mixed volume of a tuple $\AA=(A_1,A_2,\ldots, A_n)$:
\begin{equation}\label{recursive}
\prod \AA = \sum_{\a\in \ZZ} l_{A_1}(\a) \prod_{j=2}^n A_j^{\a}.
\end{equation}

\begin{prf-1}
Consider the set $\Lambda$ of all primitive covectors $\b\in \ZZ$ such that $V_\a:=\prod_{j\geq 2}^n A_j^{\beta}$ is positive. We start with the proof of the following proposition. At least one of the following conditions holds:
\begin{enumerate}
\item there exists $\b\in \Lambda$ such that $l_{A_1}(\b)>\b(\kk_0), \b(\kk_1)$.
\item there exist $\b \in \Lambda$, $\kk_2\in A_1$, and $s\in \{0,1\}$ such that $\b(\kk_{s})>\b(\kk_2)>\b(\kk_{1-s})$
\item there exist $\b_1,\b_2\in \Lambda$ and $s\in \{0,1\}$ such that $\kk_{s}\in A_1^{\b_j}$ for $j=1,2$ and $\dim A_1^{\b_j}\geq 1$ for at least one $j$.
\end{enumerate}
Assume that conditions 1,2 do not hold and prove condition~3. The linear span of $\Lambda$ has dimension~$n$. Otherwise, due to equation~(\ref{recursive}), we would have $I\cdot \prod_{j\geq 2}^n A_j=0$ for a lattice segment $I$ that is orthogonal to any $\beta\in \Lambda$. This would mean that some subtuple of $(A_2,A_3,\ldots, A_n)$ has codimension not exceeding $0$, which contradicts the conditions of Lemma~\ref{lem-main}. We also have $\sum_{\beta\in \Lambda} (\prod_{j\geq 2}^n A_j^{\beta}) \beta =0$, therefore, $|\Lambda|\geq n+1\geq 3$ (we assume that $n>1$, since the case $n=1$ is covered by Theorem~\ref{Thm_equal_polytops}, see the example after its statement). Note that, for each $\beta\in \Lambda$, at least one of $\kk_0,\kk_1$ is not contained in $A_1^{\beta}$. Indeed, if $\{\kk_0, \kk_1\}\subset A_1^{\b}$, then $[\kk_0,\kk_1]\in A_1^{\b}$, which contradicts the conditions of Lemma~\ref{lem-main}. It follows that there exists $s\in\{0,1\}$ and two different $\b_1,\b_2\in \Lambda$ such that $l_{\b_j}(A_1)>\b_j(\kk_{1-s})$, $j=1,2$. Since condition 1 does not hold, we have $\kk_s\in A_1^{\b_j}$, $j=1,2$. Conditions of Lemma~\ref{lem-main} imply that $\dim A_1>1$. Therefore there exists $\kk_2\in A_1\setminus\{\kk_0, \kk_1\}$. Since $\sum_{\beta\in \Lambda} (\prod_{j\geq 2}^n A_j^{\beta}) \beta(\kk_2-\kk_{1-s}) =0$ and the linear span of $\Lambda$ has dimension $n$, there exists $\b\in \Lambda$ such that $\b(\kk_2-\kk_{1-s})>0$ ($\b$ may either belong or not belong to $\{\b_1, \b_2\}$). Since conditions 1, 2 do not hold, we have $\b(\kk_2)=\b(\kk_{s})=l_{A_1}(\b)$, that is, $\kk_2, \kk_{s}\in A_1^\b$. It follows that $\dim A_1^{\b}\geq 1$, which completes the proof of condition 3.

We choose $s\in \{1,2\}$ with respect to conditions~2,3 (in the case when condition~1 holds, we choose it arbitrarily). We further denote $\kk_{1-s}=\kk$, assume w.l.o.g. that $\kk_{s}=0$, and prove that group $G:=G_0=G_{\kk_s}$ coincides with $S_D$. Lemma~\ref{transitivity} guarantees that the monodromy group $G$ is transitive. We prove that the monodromy group $G$ is generated by several transpositions and several permutations $a_1,a_2,\ldots, a_t\in S_D$ such that at least one of the following conditions holds:
\begin{enumerate}
\item We have $t>1$, and $a_1,a_2,\ldots, a_t$ are disjoint permutations.
\item We have $t=1$ and there is at least one $x$ such that $a_1(x)=x$. 
\end{enumerate}
In both cases, Lemma~\ref{lem-main} follows from Lemma~\ref{Permutations}.

Consider $f_{1,\lambda}=f_1+ \lambda \xx^{\kk}$. Let $\p_\lambda\colon X\to \C$ be given by $\p_\lambda(\xx)=f_{1,\lambda}(\xx)$. For generic $\lambda$, there is a correctly defined monodromy group $G_{\lambda}$ of $\p_\lambda$ that acts on $F_{\lambda}=\p_{\lambda}^{-1}(0)$. The group $G_\lambda$ is evidently isomorphic to $G$, and, instead of $G$, we will study $G_\lambda$ for large values $|\lambda|$ using the following description of the bifurcation set of the covering~$\pi_{\lambda}$ .

Consider a toric compactification $M\supset (\C\setminus\{0\})^n$ that corresponds to the set of polytopes $\conv A_1,$ $\conv A_2,\ldots, \conv A_n$. For generic $(f_2,\ldots, f_n)\in \C^{A_2}\oplus \ldots \oplus \C^{A_n}$, the closure $\overline{X}\subset M$ of the set $X$ is transversal to $M':=M\setminus (\C\setminus\{0\})^n$. There is a finite set of covectors $\a\in \Lambda$ such that $0\in A_1^\alpha$. Denote them by $\{\a_1,\a_2,\ldots, \a_g\}=\Lambda' \subset \Lambda$. The curve $\overline{X}$ intersects the stratum $M_j:=M_{\a_j}\subset M$ at $V_{j}:=V_{\alpha_j}$ points $x_{j,1}, x_{j,2}, \ldots, x_{j,V_j}\in M_j$. For each $j$ such that $\dim A_1^{\alpha_j}>0$, denote $c_{j,k}=f_{1,\lambda}(x_{j,k})=f_1(x_{j,k})\in \C$, $k=1,2,\ldots, V_j$. For generic $(f_1,f_2,\ldots, f_n)$, the values $c_{j,k}$ do not coincide and are not equal to $0$. The set of bifurcation points $B\subset \C$ of the covering $\p_{\lambda}$ can be represented as $B=B_1\sqcup B_2$, where $B_1$ contains the values $c_{j,k}$ for some $j,k$ and $B_2=\{c\in \C\mid c \mbox{ is a critical value of } {f_{1,\lambda}|}_X\}$. In particular, $B_1$ includes $B_1':=\{c_{j,k}\mid h_j>1\}$, where $h_j:=\min_s(l_{\a_j}(A_s)-l_{\a_j}(A_s\setminus A_s^{\a_j}))$ is the number of points of $\p_{\l}^{-1}(c)$ that tend to $x_{j,k}$ as $c$ tends to $c_{j,k}$. If there exists $j$ such that $A_1^{\a_j}=\{0\}$ and $h_j>1$, then $B_1=\{f_1(0)\}\sqcup B_1'$ (for generic $\{f_j\}$, we have $f_1(0)\neq 0$ and $f_1(0)\neq c_{j,k}$ for any $j,k$). Otherwise, $B_1=B_1'$. Now we describe a set of generators of the group $G_{\lambda}$.

Consider a set of disjoint loops $\{s_c\colon [0,1]\to \C\mid c\in B\}$ at point $0\in \C$, where each loop $s_{c}$ goes around the value $c\in B$ and does not link any other bifurcation values of $\p_\lambda$. The monodromy group $G_{\lambda}$ is generated by transformations of loops $s_{c}$, $c\in B$. The monodromy transformation of a loop $s_{c}$, where $c\in B_2$, is a transposition. Now we show that the transformations of the loops $s_c$, $c\in B_1$, are disjoint permutations.

Consider a small enough tubular neighborhood $U\subset M$ of $M'$. For sufficiently large $|\lambda|$, we have $\{f_{1,\lambda}=0\}\subset U$. We can show that the transformation of $s_{c_{j,k}}$ acts only on the set $U_{j,k}\cap F_{\lambda}$, where $U_{j,k}\subset U$ is a neighborhood of~$x_{j,k}$ and the cardinality of the set is $H_j:=|U_{j,k}\cap F_{\lambda}|= l_{A_1}(\a_j)-\a_{j}(\kk)=-\a_{j}(\kk)$. In fact, for sufficiently large $|\lambda|$, we have $\{f_{1,\lambda}=c_0\}\subset U$ for any $c_0\in \cup_{c\in B_1} s_c([0,1])$. This implies that the monodromy transformation of $s_{c}$ permutes some $h_j$ of $H_{j}$ points of $F_{\lambda}\cap U_{j,k}$ by a cycle of length $h_j$ and does not move other points of $F_{\lambda}$. If $|B_1|>1$, Lemma~\ref{lem-main} immediately follows now from Lemma~\ref{Permutations}.

Now assume that $|B_1|=1$. Denote by $c$ the only point of $B_1$. We will show that the transformation $a\colon F_{\l}\to F_{\l}$ of the loop~$s$ going around $c$ has a fixed point. Denote $\a_{j_1},\a_{j_2},\ldots, \a_{j_t}$ all covectors $\g\in\Lambda$ such that $\{f_1^{\g}-c=f_2^{\g}=\ldots=f_n^{\gamma}=0\}$ has a solution in $(\C^*)^n$, where $f_j^{\g}$ is the part of the Laurent polynomial $f_j$ formed by monomials corresponding to integer points of $A_j^\g$ (note that $t=1$, if $c\neq f_1(0)$). The number of points that are not fixed under $a$ can be estimated from above by
\begin{multline}
\sum_{s=1}^t h_{j_s}V_{j_s}\leq \sum_{s=1}^t H_{j_s}V_{j_s} = \sum_{s=1}^t (l_{\a_{j_s}}(A_1)-\a_{j_s}(\kk))\prod_{p=2}^n A_p^{\a_{j_s}}\leq\sum_{\b\in \Lambda} (l_{\b}(A_1)-\b(\kk))\prod_{j=2}^n A_j^{\b}=\\
=\sum_{\b\in \Lambda} l_{\b}(A_1)\prod_{j=2}^n A_j^{\b} - \sum_{\b\in \Lambda} \b(\kk)\prod_{j=2}^n A_j^{\b}= \sum_{\b\in \Lambda} l_{\b}(A_1)\prod_{j=2}^n A_j^{\b} = \prod_{j=1}^n A_j = |F_{\l}|.\\
\end{multline}
In the case when condition 2 holds, the first inequality is strict. In the cases when one of the conditions 1,3 holds, the second inequality is strict. In any case, we have $\sum_{s=1}^t h_{j_s}V_{j_s}<|F_{\l}|$, therefore, there is a fixed point of $a$ and Lemma~\ref{lem-main} follows from Lemma~\ref{Permutations}.
\end{prf-1}

%{\noindent(A. Esterov) National Research University Higher School of Economics. \newline Faculty of Mathematics NRU HSE, 7 Vavilova 117312 Moscow, Russia.}

%{\noindent(G. Gusev) Moscow Institute of Physics and Technology (State University). \newline Department of Innovations and High Technology, 9 Institutskii per. \newline 141700 Dolgoprudny, Moscow Region, Russia.}

\end{document}